\documentclass[final]{siamltex}
\usepackage[dvips]{graphicx}
\usepackage{amsmath}
\usepackage{latexsym}
\usepackage{amssymb}
\usepackage{algorithm,algorithmic}
\usepackage{color}
 \usepackage{mathptmx}

\title{Numerical  solutions to large-scale differential Lyapunov matrix equations}
\author{ M. Hached \thanks{Laboratoire P.  Painlev\'e UMR 8524, UFR de Math\'{e}matiques, Universit\'e des Sciences et Technologies de Lille, IUT A, Rue de la Recherche, BP 179, 59653 Villeneuve d'Ascq Cedex, France; {\tt email:  mustapha.hached@univ-lille1.fr}} \and K. Jbilou\thanks{L.M.P.A, Universit\'e du Littoral C\^ote d'Opale,
50 rue F. Buisson BP. 699, F-62228 Calais Cedex, France. {\tt E-mail: ;
jbilou@univ-littoral.fr}.} }
\date{ }
\newtheorem{Remark}{Remark}

\newcommand{\K}{\mathbb{K}}

\begin {document}
\maketitle
\begin{abstract}
In the present paper, we consider large-scale  differential Lyapunov matrix  equations having a low rank constant term. We present two new approaches for the numerical resolution of such differential matrix equations. The first approach is based on the integral expression of the exact solution and an approximation method for  the computation of the exponential of a matrix times a block of vectors.   In the second approach, we first  project the initial problem onto a block (or extended block)  Krylov subspace and get a low-dimensional differential Lyapunov matrix  equation. The latter differential matrix problem is then solved by the Backward Differentiation Formula method (BDF) and the obtained solution is used to build the low rank approximate solution of the original problem. The process being repeated until some prescribed accuracy is achieved. We give some new  theoretical results and present some numerical experiments.
\end{abstract}
\pagestyle{myheadings}
\thispagestyle{plain}
\markboth{M. Hached and   K. Jbilou}{Low rank approximate solutions ...}
\section{Introduction}
In the present paper, we consider the  differential Lyapunov 
 matrix  equation (DLE in short) of the form
\begin{equation}\label{lyap1}
\left\{
\begin{array}{l}
\dot X(t)=A(t)\,X(t)+X(t)\,A^T(t)+B(t)B(t)^T;\; (DLE) \\
 \;X(t_0)=X_0,\; \; t \in [t_0, \, T_f],
\end{array}
\right.
\end{equation}

\noindent where the matrix $ A(t) \in \mathbb{R}^{n \times n} $ is assumed to be nonsingular  and  $ B(t) \in
\mathbb{R}^{n \times s} $  is  a full rank matrix, with $ s \ll n $.  The initial condition $X_0$ is assumed to be a symmetric and positive low-rank given matrix. \\
Differential Lyapunov  equations play a fundamental role in many areas such as control,
filter design theory, model reduction problems, differential
equations and robust control problems \cite{abou03,corless}.
For those  applications, the matrix  $A$ is generally sparse 
and very large. For such problems, only  a few attempts have been
made to solve (\ref{lyap1}). \\

Let us first recall the following theoretical result which gives an expression of the exact solution of \eqref{lyap1}.
\begin{theorem}\label{theo1}
\cite{abou03}
The unique solution of the  general Lyapunov differential equation
\begin{equation}
\displaystyle {\dot X}(t)=A(t)\,X+X\,A(t)^T+M(t);\;\; X(t_0)=X_0
\end{equation}
is defined by
\begin{equation}\label{solexacte1}
X(t) = \Phi_A(t,t_0)X_0\Phi^T_A(t,t_0)+\int_{t_0}^t \Phi_A(t,\tau)M(\tau)\Phi^T_A(t,\tau)d\tau.
\end{equation}
where the transition matrix $\Phi_A(t,t_0)$ is the unique solution to the problem
$$\displaystyle {\dot \Phi}_A(t,t_0)=A(t) \Phi_A(t,t_0),\;\; \Phi_A(t_0,t_0)=I.$$
Futhermore, if $A$ is assumed to be a constant matrix, then we have
\begin{equation}\label{solexacte2}
X(t)=e^{(t-t_0)A}X_0e^{(t-t_0)A^T}+\int_{t_0}^t e^{(t-\tau)A}M(\tau)e^{(t-\tau)A^T}d\tau.
\end{equation}
\end{theorem}
We notice that the problem \eqref{lyap1} is equivalent to the linear ordinary differential equation
\begin{equation}\label{kron}
\left\{
\begin{array}{c c l}
\dot{x}(t)& =& \mathcal{A}(t)x(t)+b(t) \\
x_0 & = & vec(X_0)
\end{array}
\right.
\end{equation}
where $\mathcal{A}= I \otimes A(t) + A(t) \otimes I$, $x(t)=vec(X(t))$ and $b(t) = vec(B(t)B(t)^T)$, where $vec(Z)$ is the long  vector obtained by stacking the columns of the matrix $Z$. For moderate size problems, it is then possible to use an integration method to solve  \eqref{kron}. However, this approach is not adapted to large problems. In the present paper, we will consider projection methods onto extended block Krylov (or block Krylov if $A$ is not invertible) subspaces  associated to the pair $(A,B)$. These subspaces are defined as follows
$$\K_m(A,B)={\rm range}(B,AB,\ldots,A^{m-1}B)$$
for block Krylov subspaces, or
$${\cal K}_m(A,B)={\rm range}(A^{-m},\ldots,A^{-1}B,B,AB,\ldots,A^{m-1}B)$$
for extended block Krylov subspaces. Notice that  the extended Krylov subspace ${\cal K}_k(A,B)$ is a sum of two block Krylov subspaces
$$ {\cal K}_m(A,B)=\K_{m}(A,B)\, + \, \K_m(A^{-1},A^{-1}B). $$
To compute an orthonormal basis $\{V_1,\ldots,V_m\}$,  where $V_i$ is of dimension $n\times s$ for the block Krylov and $n \times 2s$ in the extended  block Krylov case, two  algorithms have been defined: the first one is the well known block Arnoldi algorithm and the second one is the extended block Arnoldi algorithm \cite{druskin98,simoncini1}. These algorithms also generate block Hessenberg matrices ${\bar {\cal T}_m} = {\cal V}_{m+1}^T\,A\,{\cal V}_m $ satisfying the following algebraic relations
\begin{eqnarray}
A\,{\cal V}_m & = & {\cal V}_{m+1}\,{\bar {\cal T}}_m, \\
              & = & {\cal V}_m\,{\cal T}_m + V_{m+1}\,T_{m+1,m}\,E_m^T,
\end{eqnarray}
where ${\cal T}_m = {\bar {\cal T}_m} (1:d,:)={\cal V}_m^T A {\cal V}_m$ and 
where $ T_{i,j} $ is the  $ (i,j) $ block of $  {\bar {\cal T}_m} $ of size $ d \times d $, and $ E_m = [ O_{d \times (m-1)d}, I_{d} ]^T $ is the matrix of the last $ d $ columns of the $ md \times md$ identity matrix $ I_{md} $ with $d=s$ for the block Arnoldi and $d=2s$ for the extended block Arnoldi. \\
When the matrix $A$ is nonsingular and when the computation of $W=A^{-1}V$ is not difficult (which is the case for sparse and structured matrices), the use of the extended block Arnoldi is to be preferred. \\  
The paper is organized as follows: In Section 2, we present a first approach based on the approximation of the exponential of a matrix times a block using a Krylov projection method. We give some theoretical results such as an upper bound for the norm of the error and an expression of the exact residual.  A second approach,presented in Section 3, for  which the initial differential Lyapunov matrix equation is projected onto a block (or extended block) Krylov subspace. Then, the obtained low dimensional differential Lyapunov equation is solved by using the well known Backward Differentiation Formula (BDF). In Section 4, an application to balanced truncation method for large scale linear-time varying dynamical systems is presented. The last section is devoted to some numerical experiments.


\section{The first  approach: using an approximation of the matrix exponential}

In this section, we give a new  approach for computing approximate solutions to large differential equations \eqref{lyap1}. 
The expression of the exact solution as
\begin{equation}\label{solexacte3}
X(t)=e^{(t-t_0)A}X_0e^{(t-t_0)A^T}+\int_{t_0}^t  e^{(t-\tau)A}\, B B^Te^{(t-\tau)A^T}\, d\tau,
\end{equation}
suggests the idea of computing $X(t)$  by approximating the factor $ e^{(t-\tau)A}B$ and then  using a quadrature method to compute the desired approximate solution. \\ As computing the exponential of a small matrix is straightforward , this is not the case for large scale problems, as   $e^{(t-\tau)A}$ could be dense even though $A$ is sparse. However, in our problem, the computation of $e^{(t-\tau)A}$ is not needed as we will rather consider the product $e^{(t-\tau)A}\, B$, for which approximations via projection methods onto block or extended block Krylov subspaces are well suited. 
\\
Krylov subspace projection methods generate a sequence of  nested subspaces (Krylov or extended Krylov subspaces). 
Let $\mathcal{V}_m=[V_1,\ldots,V_m]$ be the orthogonal  matrix whose  columns form an orthonormal  basis of the subspace $\mathit{K}_m$,  
Following \cite{saad1,saad2,vorst1}, an approximation to $Z=e^{(t-\tau)A}\, B$ can be obtained as
\begin{equation}
\label{exp1}
Z_m(t) = \mathcal{V}_m e^{(t-\tau)\mathcal{T}_m}\, \mathcal{V}_m^T B
\end{equation}
where $\mathcal{T}_m=\mathcal{V}^T_m A \mathcal{V}_m$. 
Therefore, the term appearing in the integral expression \eqref{solexacte3} can be approximated as 
\begin{equation}
\label{exp2}
e^{(t-\tau)A}BB^Te^{(t-\tau)A^T} \approx Z_m(t) Z_m(t)^T.
\end{equation}
If for simplicity, we assume $X_{t_0}=0$, an approximation to the solution of the differential Lyapunov equation \eqref{solexacte3} can be expressed as 
\begin{equation}
\label{exp3}
X_m(t) = \mathcal{V}_m G_m(t) {\mathcal{V}_m}^T,
\end{equation}
where 
\begin{equation}
\label{gm}
G_m(t) = \displaystyle  \int_{t_0}^t  {\widetilde G}_m(\tau) {\widetilde G}_m^T(\tau) d\tau,
\end{equation}
 and $ {\widetilde G}_m(\tau)= e^{(t-\tau)\mathcal{T}_m}B_m$. \\
The next result shows that the matrix function $G_m$ is the solution of a low-order differential Lyapunov matrix equation.
\begin{theorem}
Let $G_m(t)$ be the matrix function  defined by \eqref{gm}, then it satisfies the following low-order differential Lyapunov matrix equation
\begin{equation}\label{low2}
{\dot G}_m(t) = \mathcal{T}_m G_m(t) + G_m(t){\mathcal{T}_m}^T  +B_mB_m^T,\; t \in [t_0,\,T_f]
\end{equation}
\end{theorem}
\begin{proof}
The proof can be easily derived from the expression \eqref{gm} and the result of Theorem \ref{theo1}.
\end{proof}

\medskip
As a consequence, intruducing  the residual $ R_m(t) = \displaystyle {\dot X}_m(t)-A\,X_m-X_m\,A^T- BB^T $ associated to the approximation $X_m$,  we have the following relation
\begin{eqnarray*}
{\cal V}^T_m R_m(t) {\cal V}_m &= & {\cal V}^T_m ({\dot X} -AX_m(t)-X_m(t)A^T-BB^T) {\cal V}_m\\
& = &  {\dot G}_m(t) -\mathcal{T}_m G_m(t) - G_m(t){\mathcal{T}_m}^T  -B_mB_m^T\\
& = & 0,
\end{eqnarray*}
which shows that the residual satisfies a Petrov-Galerkin condition.\\

As mentioned earlier,  once  $ {\widetilde G}_m(\tau)$ is computed, we use a quadrature method to approximate the integral \eqref{gm} in order to approximate $G_m(t)$.\\
We now briefly discuss  some practical aspects of the computation of $e^{(t-\tau)\mathcal{T}_m} B_m$ where $B_m={\cal V}_m^TB$, when $m$ is small and $\mathcal{T}_m$ is a an upper block Hessenberg matrix.\\

\noindent In the last decade, many approximation techniques such as the use of partial fraction expansions or Pad\'e approximation have been proposed, see for example  \cite{gallopoulos,saad2}. However, it was remarked  that a good way for evaluating the exponential of matrix times by a vector by using rational approximation to the exponential function. One of the main advantages of  rational approximations as compared to polynomial approximations is  the better stability of their  integration schemes. 
Let us consider the rational function 
$$F(z)= a_0 + \displaystyle \sum_{i=1}^p \frac{a_i}{z-\theta_i},$$
where the $\theta_i$'s are the poles of the rational function $F$. Then, the approximation to ${\widetilde G}_m(\tau)=e^{(t-\tau)\mathcal{T}_m} $ is given by
\begin{equation}
\label{ratapp1}
{\widetilde G}_m(\tau) \approx a_0B_m +  \displaystyle \sum_{i=1}^p  a_i [{(t-\tau)\mathcal{T}_m}- \theta_i I]^{-1} \, B_m.
\end{equation}
One of the possible choices for the rational function $F$ is  based on Chebychev approximation  of the function $e^x$ on $[0,\, \infty[$, see \cite{saad2}.  We  notice that for small values of $m$,  one can also directly compute the matrix exponential  $e^{(t-\tau){\mathcal T}_m}$ by using the well-known 'scaling and squaring method for the matrix exponential' method, \cite{higham09}. This method was  associated to  a Pad\'e approximation and is implemented in the {\tt expm} Matlab routine. \\

\noindent From now on, we assume that the basis formed by the orthonormal  columns of ${\cal V}_m$  is obtained by applying  the   block Arnoldi or the extended block Arnoldi algorithm  to the pair $(A,B)$. \\
The computation of $ X_m(t) $
(and of $ R_m(t) $) becomes expensive as $ m $ increases. So, in
order to stop the iterations, one has to test if $ \parallel R_m
\parallel < \epsilon $ without having to compute extra products
involving the matrix $ A $. The next result shows how to compute
the residual norm of $ R_m(t) $ without forming the approximation $
X_m(t) $ which is computed in a factored form only when convergence
is achieved. 
\medskip
\begin{theorem} \label{t2}
Let $ X_m(t) = {\cal V}_mG_m(t){\cal V}_m^T $ be the approximation obtained at step $ m $ by the  block (or extended block) 
Arnoldi  method. Then the residual $ R_m(t) $ satisfies

\begin{equation}
\label{result2}
\parallel R_m(t) \parallel = \parallel T_{m+1,m} \bar G_m(t) \parallel,
\end{equation}
where $ \bar G_m $ is the $ d \times md $  matrix corresponding to the last $ d $ rows of $ G_m $ where $d=s$ when using the block Arnoldi and $d=2s$ for the extended block Arnoldi.
\end{theorem} 
\medskip
\begin{proof}
The proof of this theorem comes directly from \eqref{exp3} and the fact that $G_m$ solves the low dimensional problem \eqref{low2}.
\end{proof}

\medskip
\noindent The result of Theorem \ref{t2} is very  important in practice, as it allows us to
stop the iterations when convergence is achieved without computing the approximate solution $X_m(t)$. \\ The following result  shows  that the approximation $X_m $ is an exact solution of a perturbed  differential  Lyapunov  equation. \\

\begin{theorem}
Let $X_m(t)$ be the approximate solution given by \eqref{exp3}. Then we have 
\begin{equation}
\label{pertu}
\displaystyle {\dot X}_m(t)=(A-F_m)\,X_m+X_m\,(A-F_m)^T+BB^T.
\end{equation}
where $ F_m = V_m\,T_{m+1,m}^T\,V_{m+1}^T $.
\end{theorem}

\medskip
\begin{proof}
The proof is easily obtained  from \eqref{low2} and  the expression \eqref{exp3} of the  approximate solution $X_m(t)$.\\
\end{proof}

\medskip
\begin{Remark}
The solution
$X_m(t)$   can be given as a product of two low rank matrices.
Consider the eigen-decomposition of the symmetric and positive matrix 
$md \times md$ $G_m(t)=U\, D\, U^T$ where $D$ is the
diagonal matrix of the  eigenvalues of $G_m(t)$ sorted in
decreasing order and $d=s$ for the block Arnoldi or $d=2s$ for the extended block Arnoldi. Let $U_l$ be the $md \times l$ matrix  of   the first $l$ columns of  $U$ 
corresponding to the $l$ eigenvalues of magnitude greater than
 some tolerance $dtol$. We obtain the
truncated eigen-decomposition  $G_m(t) \approx U_l\, D_l\,
U_l^T$ where $D_l = {\rm diag}[\lambda_1, \ldots, \lambda_l]$.
 Setting ${\widetilde Z}_m(t)={\cal V}_m \, U_l\, D_l^{1/2}$, it
follows that
\begin{equation}
\label{approx}
X_m(t) \approx {\widetilde Z}_m(t) {\widetilde Z}_m(t)^T.
\end{equation}
Therefore, one has to compute and to store only the matrix ${\widetilde Z}_m(t)$ which is usually the required factor  in some control problems such as in the balanced truncation method for model reduction in large scale dynamical systems. This possibility is very important for storage limitations in the  large scale problems.
\end{Remark}

\noindent The next result states that the error matrix $X(t)-X_m(t)$ satisfies a differential Lyapunov matrix equation.
\medskip
\begin{theorem}\label{err1}
Let $X(t)$ be the exact solution of \eqref{lyap1} and let $X_m(t)$ be the approximate solution obtained at step $m$. The error $E_m(t)=X(t)-X_m(t)$  satisfies the following equation
\begin{equation}\label{pertu2}
\displaystyle {\dot E}_m(t)=AE_m(t)+E_m(t)A^T-R_m(t),
\end{equation}
and
\begin{equation}
\label{error3}
E_m(t)=e^{(t-t_0)A}E_{m,0}e^{(t-t_0)A^T}+\int_{t_0}^t e^{(t-\tau)A}R_m(\tau)e^{(t-\tau)A^T}d\tau,\; t \in [t_0,\, T_f].
\end{equation}
where  $E_{m,0}=E_m(0)$.
\end{theorem}
\begin{proof}
The result is easily obtained by subtracting the residual equation from the initial differential Lyapunov equation  \eqref{lyap1}.
\end{proof}
\medskip

\medskip
\noindent Next, we give  an upper bound for the norm of the error in the case where $A$ is a stable matrix.\\

\begin{theorem}
\label{Theoerr2}
Assume that $A$ is a  stable matrix and $X(t_0)=X_m(t_0)$. Then we have the following upper bound
\begin{equation}
\label{upperbound}
\parallel E_m(t) \parallel  \le   \displaystyle  \parallel  T_{m+1,m} \parallel \, \parallel \bar G_m \parallel_{\infty}  \frac{e^{2(t-t_0)\mu_2(A)}-1}{2 \mu_2(A)},\\
\end{equation}
where $\mu_2(A)=\displaystyle \frac{1}{2} \lambda_{max}(A+A^T)<0$  is the 2-logarithmic norm   and $\parallel \bar G_m \parallel_{\infty}  =\displaystyle \max_{\tau \in [t_0,\, t]} \parallel \bar G_m(\tau) \parallel$. The matrix  $ \bar G_m $ is the $ d \times md $  matrix corresponding to the last $ d $ rows of $ G_m $ where $d=s$ when using the block Arnoldi and $d=2s$ for the extended block Arnoldi.
\end{theorem}
\medskip
\begin{proof}
We first remind that if $A$ is a stable matrix, then the logarithmic norm provides the following bound  $\parallel e^{tA} \parallel \le e^{\mu_2(A)t}$. Therefore, using the expression \eqref{error3}, we  obtain the following relation 
\begin{equation*}
\parallel E_m(t) \parallel  \le   \displaystyle \int_{t_0}^t \parallel e^{(t-\tau)A} \parallel^2 \, \parallel R_m(\tau) \parallel d \tau.
\end{equation*}
Therefore, using \eqref{result2} and the fact that $\parallel e^{(t-\tau)A} \parallel \le e^{(t-\tau) \mu_2(A)}$, we get 
\begin{eqnarray*}
\parallel E_m(t) \parallel &  \le  &  \parallel  T_{m+1,m}  \bar G_m \parallel_{\infty}  \displaystyle \int_{t_0}^t e^{2(t-\tau) \mu_2(A)} d\tau\\
 & \le & \parallel  T_{m+1,m} \parallel \, \parallel \bar G_m \parallel_{\infty} e^{2t\mu_2(A)} \displaystyle \int_{t_0}^t e^{-2\tau \mu_2(A)} d\tau\\
 & \le & \parallel  T_{m+1,m}  \parallel \parallel \bar G_m \parallel_{\infty} e^{2t\mu_2(A)} \, \times  \frac{e^{-2\mu_2(A)t}-e^{-2\mu_2(A)t_0}}{-2 \mu_2(A)}\\
 & = & \displaystyle  \parallel  T_{m+1,m} \parallel \, \parallel \bar G_m \parallel_{\infty} \frac{e^{2(t-t_0)\mu_2(A)}-1}{2 \mu_2(A)},
\end{eqnarray*}
which gives the desired result.

\end{proof}

\medskip
\noindent Notice that if $\parallel T_{m+1,m} \parallel $ is close to zero, which is the case when $m$ is close to the degree of the minimal polynomial of $A$ for $B$, then Theorem \ref{Theoerr2} shows that the error $E_m(t)$ tends to zero.\\

\noindent Next, we give another  error bound for the norm of the error for every matrix  $A$.
\medskip
\begin{theorem}
\label{Theoer3}
Let $X(t)$ be the exact solution to \eqref{lyap1} and let $X_m(t)$ be the approximate solution obtained at step $m$.  Then we have
\begin{eqnarray*}
\label{errX}
\Vert X(t)- X_m(t) \Vert & \le & \,e^{t\mu_2(A)} (\Vert B \Vert +  \Vert B_m \Vert )\int_{t_0}^t  e^{-\tau \mu_2(A)} \Vert e^{(t-\tau)A}B-{\cal V}_m e^{(t-\tau) \mathcal{T}_m} B_m\Vert d \tau
\end{eqnarray*}
where $\mu_2(A)= \lambda_{max}((A+A^T)/2)$, $Z(\tau) = e^{(t-\tau)A} B$ and $Z_m(\tau)={\cal V}_m e^{(t-\tau)\mathcal{T}_m} B_m$ with $B_m={\cal V}_m^T B$.
\end{theorem}
\medskip
\begin{proof}
From the expressions of $X(t)$ and $X_m(t)$, we have 
\begin{eqnarray*}
\Vert X(t) - X_m(t) \Vert & = & \left  \Vert \int_{t_{0}}^t (Z(\tau)Z(\tau)^T -Z_m(\tau)Z_m(\tau)^T) d \tau  \right \Vert\\
&= & \left  \Vert \int_{t_{0}}^t [Z(\tau)(Z(\tau)-Z_m(\tau))^T +(Z(\tau)-Z_m(\tau))Z^T_m(\tau)] d \tau \right \Vert\\
& \le & \int_{t_{0}}^t  (\Vert Z(\tau) \Vert + \Vert Z_m(\tau) \Vert) \Vert  Z(\tau) -Z_m(\tau) \Vert  d \tau,\\
\end{eqnarray*}
Therefore, using the fact that $\mu_2(\mathcal{T}_m)=\lambda_{max}((\mathcal{T}_m+\mathcal{T}_m^T)/2) \le  \lambda_{max}((A+A^T)/2)=\mu_2(A)$, where $\mathcal{T}_m= \mathcal{V}_m^T A \mathcal{V}_m$, it follows that 
\begin{eqnarray*}
\Vert X(t) - X_m(t) \Vert  & \le &   e^{t\mu_2(A)} (\Vert B \Vert +  \Vert B_m \Vert ) \int_{t_0}^t  e^{-\tau \mu_2(A)}\Vert Z(\tau) -  Z_m(\tau) \Vert d \tau\\
& \le & \,e^{t\mu_2(A)} (\Vert B \Vert +  \Vert B_m \Vert )\int_{t_0}^t  e^{-\tau \mu_2(A)} \Vert e^{(t-\tau)A}B-{\cal V}_m e^{(t-\tau) \mathcal{T}_m} B_m\Vert d \tau,\\
\end{eqnarray*}
\end{proof}

\medskip
\noindent When using a block Krylov subspace method such as the block Arnoldi method, then one can generalize to the block case the results already stated  in many papers; see \cite{druskin98,gallopoulos,higham05,saad2}. In particular, we can easily generalize the result given in \cite{saad2} for the case $s=1$ to the case $s >1$. In this case, we have the following upper bound.
\begin{equation}\label{y1}
\Vert e^{A}B-{\cal V}_m e^{\mathcal{T}_m} B_m \Vert \le 2 \parallel B \parallel  \; \displaystyle \frac{\rho^m e^{\rho}}{m!},
\end{equation}
where $\rho= \Vert A \Vert$
\medskip 

\noindent The rupper bound \eqref{y1} could be used in Theorem \ref{Theoer3} to obtain a new upper bound for the norm of the error. In that case, we obtain the following upper bound
\begin{equation}
\label{er1}
\Vert X(t) - X_m(t) \Vert \le 2 \parallel B \parallel \displaystyle \frac{\rho^m}{m!} \,e^{t(\mu_2(A)+\rho)} (\Vert B \Vert +  \Vert B_m \Vert )\int_{t_0}^t  e^{-\tau (\mu_2(A)+\rho)} (t-\tau)^m d \tau,\\
\end{equation}


We summarize the steps of our proposed first approach (using  the extended block Arnoldi) in the following algorithm
\begin{algorithm}[h!]
\caption{The  extended block Arnoldi (EBA-exp) method for DLE's}\label{algo_EBA_exp}
\begin{itemize}
\item Input $X_0=X(t_0)$, a tolerance $tol>0$, an integer $m_{max}$.
\item  For $ m = 1,\ldots,m_{max} $
\begin{itemize}
\item  Apply the extended block Arnoldi algorithm to compute an orthonormal basis ${\mathcal V}_m=[V_1,...,V_m]$ of ${\mathcal K}_m(A,B)=Range[B,A^{-1}B,...,A^{-m}B,A^{m-1}B]$ and the upper block Hessenberg matrix ${\mathcal T}_m$.
\item Set ${B}_m={\mathcal V}_m^TB$ and compute $ {\widetilde G}_m(\tau)= e^{(t-\tau)\mathcal{T}_m}B_m$ using the matlab function {\tt expm}.
\item Use a quadrature method to compute the integral \eqref{gm} and get an approximation of $G_m(t)$ for each $t \in [t_0,\, T_f]$.
\item If $\parallel R_m(t) \parallel = \parallel T_{m+1,m} \bar G_m(t) \parallel < tol$ stop and compute the approximate solution $X_m(t)$ in the factored form given by the relation \eqref{approx}.
\end{itemize}
\item End
\end{itemize}
\end{algorithm}

\section{A second approach: Projecting  and solving with BDF}
\subsection{Low-rank approximate solutions via BDF}\label{ss3.1}
In this section, we show how to obtain low rank approximate solutions to the differential Lyapunov equation \eqref{lyap1} by projecting directly the initial problem onto small block Krylov or extended block Krylov subspaces. \\
We first apply the  block Arnoldi  algorithm (or the extended block Arnoldi)  to the pair $(A,B)$ to get the matrices ${\cal V}_m$ and $ {\cal T}_m={\cal V}^T_m A {\cal V}_m $. Let $X_m(t)$ be the desired low rank approximate solution  given as 
\begin{equation}\label{approx1}
X_m(t) = {\cal V}_m Y_m(t) {\cal V}_m^T,
\end{equation}
satisfying the Petrov-Galerkin orthogonality condition
\begin{equation}
\label{galerkin}
{\cal V}_m^T R_m(t) {\cal V}_m =0,\; t \in [t_0,\; T_f],
\end{equation}
where $R_m(t)$ is the residual $ R_m(t) = \displaystyle {\dot X}_m(t)-A\,X_m(t)-X_m(t)\,A^T- BB^T $.  Then, from \eqref{approx1} and \eqref{galerkin}, we obtain the low dimensional differential Lyapunov  equation
\begin{equation}\label{lowlyap}
\displaystyle {\dot Y}_m(t)- {\cal T}_m\,Y_m(t)-Y_m(t)\,{\cal T}_m^T  - B_mB_m^T=0,
\end{equation}
with  ${\cal T}_m= {\cal V}^T_m A {\cal V}_m$ and  $ { B}_m= {\cal V}_m^T\,B $. The obtained low dimensional differential Lyapunov equation \eqref{lowlyap} is  the 
same as the one given by \eqref{low2}. For this second approach, we have to solve the latter low dimensional differential Lyapunov equation by some integration method such as the well known  Backward Differentiation  Formula (BDF).\\
Notice that we can also compute the norm of the residual without computing the approximation $X_m(t)$ which is also given, when convergence is achieved, in a factored form as in \eqref{approx}. The norm of the residual is given as
\begin{equation}
\label{result22}
\parallel R_m(t) \parallel = \parallel T_{m+1,m} \bar Y_m(t) \parallel,
\end{equation}
where $ \bar Y_m $ is the $ d \times md $  matrix corresponding to the last $ d $ rows of $ Y_m $ where $d=s$ when using the block Arnoldi and $d=2s$ for the extended block Arnoldi.

\subsection{BDF for solving the low order differential Lyapunov equation \eqref{lowlyap}}\label{projbdf}

In this subsection, we will apply the  Backward Differentiation Formula (BDF) method for solving, at each step $m$ of the block (or extended) block  Arnoldi process,  the low dimensional differential Lyapunov matrix equation \eqref{lowlyap}.  We notice that BDF is  especially used for the solution of stiff differential equations.\\  At each time $t_k$, let  $Y_{m,k}$ of the approximation of $Y_m(t_k)$, where $Y_m$ is a  solution of  (\ref{lowlyap}).  Then, the new approximation $Y_{m,k+1}$ of  $Y_m(t_{k+1})$ obtained at step $k+1$ by BDF is defined  by the implicit relation 
\begin{equation}
\label{bdf}
Y_{m,k+1} = \displaystyle \sum_{i=0}^{p-1} \alpha_i Y_{m,k-i} +h_k \beta {\mathcal F}(Y_{m,k+1}),
\end{equation} 
where $h_k=t_{k+1}-t_k$ is the step size, $\alpha_i$ and $\beta_i$ are the coefficients of the BDF method as listed  in Table \ref{tab1} and ${\mathcal F}(X)$ is  given by 
$${\mathcal F}(Y)= {\cal T}_m\,Y+Y\,{\cal T}_m^T+\,B_m\,B_m^T.$$

\begin{table}[h!!]
\begin{center}
\begin{tabular}{c|cccc} 
\hline
$p$ & $\beta$ &$\alpha_0$ & $\alpha_1$ & $\alpha_2$ \\
\hline
1 & 1 & 1 & &\\
2 & 2/3 & 4/3& -1/3 &\\
3 & 6/11 & 18/11 & -9/11 & 2/11\\
\hline
\end{tabular}
\caption{Coefficients of the $p$-step BDF method with $p \le 3$.}\label{tab1}
\end{center}
\end{table}
\noindent The approximate $Y_{m,k+1}$ solves the following matrix equation
\begin{equation*}
-Y_{m,k+1} +h_k\beta ({\cal T}_m Y_{m,k+1} + Y_{m,k+1} {\cal T}_m^T)+ B B^T + \displaystyle \sum_{i=0}^{p-1} \alpha_i Y_{m,k-i} = 0,
\end{equation*}
which can be written as the following  Lyapunov matrix equation

\begin{equation}
\label{lyapbdf}
\mathbb{T}_m\, Y_{m,k+1}  + \,Y_{m,k+1} \mathbb{T}_m^T+ \mathbb{B}_{m,k}\, \mathbb{B}_{m,k}^T   =0.
\end{equation}
We assume  that at each time $t_k$, the approximation $Y_{m,k}$ is  factorized  as a low rank product  $Y_{m,k}\approx Z_{m,k} {Z_{m,k}}^T$, where $Z_{m,k} \in \mathbb{R}^{n \times m_k}$, with $m_k \ll n$. In that case, the coefficient matrices appearing in \eqref{lyapbdf} are given by
$$\mathbb{T}_m= h_k\beta {\cal T}_m -\displaystyle \frac{1}{2}I \;  \mbox{and} \; \mathbb{B}_{m,k+1}=[\sqrt{h_k\beta} B^T, \sqrt{\alpha_0}Z_{m,k}^T,\ldots,\sqrt{\alpha_{p-1}} Z_{m,k+1-p}^T]^T.$$
The  Lyapunov matrix  equation \eqref{lyapbdf} can be solved by applying direct methods based on Schur decomposition such as the Bartels-Stewart algorithm \cite{bartels,gnv}. We notice that for large problems, many Krylov subspace type methods have been proposed to solve \eqref{lyapbdf}; \cite{Elguen,hureichel,jaim,jbilou,jbilou-Riquet,simoncini1,saad1}.\\

\begin{Remark}The main difference between Approach 1 and Approach 2 is the fact that in the first case, we compute an approximation of an integral using a quadrature formulae while in the second case, we have to solve a low dimensional differential Lyapunov equation using the BDF method. Mathematically, the two approaches are equivalent and they differ only in the way of computing numerically the low-order approximations:  $G_m$ in the first approach and $Y_m$ in the second one. 
\end{Remark}`

\noindent We summarize the steps of our proposed first approach (using  the extended block Arnoldi) in the following algorithm
\begin{algorithm}[h!]
\caption{The  extended block Arnoldi (EBA-BDF) method for DLE's}\label{algo_EBA_BDF}
\begin{itemize}
\item Input $X_0=X(t_0)$, a tolerance $tol>0$, an integer $m_{max}$.
\item  For $ m = 1,\ldots,m_{max} $
\begin{itemize}
\item  Apply the extended block Arnoldi algorithm to compute an orthonormal basis ${\mathcal V}_m=[V_1,...,V_m]$ of ${\mathcal K}_m(A,B)=Range[B,A^{-1}B,...,A^{-m}B,A^{m-1}B]$ and the upper block Hessenberg matrix ${\mathcal T}_m$.
\item Set ${B}_m={\mathcal V}_m^TB$ and use the BDF method to solve the low dimensional differential Lyapunov equation
$$\displaystyle {\dot Y}_m(t)- {\cal T}_m\,Y_m(t)-Y_m(t)\,{\cal T}_m^T  - B_mB_m^T=0,\; t \in [t_0,\,T_f]$$
\item If $\parallel R_m(t) \parallel = \parallel T_{m+1,m} \bar Y_m(t) \parallel < tol$ stop and compute the approximate solution $X_m(t)$ in the factored form given by the relation \eqref{approx}.
\end{itemize}
\item End
\end{itemize}
\end{algorithm}

\section{Application: Balanced truncation for linear time-varying dynamical  systems}
In this section, we assume that the coefficient matrices $A$ and $B$ are time-dependent. It is the case for example when we are dealing with Multi-Input Multi-Output (MIMO)  linear-time varying (LTV) dynamical systems
\begin{equation}
\label{ltv1}
\left\{ \begin{array}{lcl}
 \dot {x}(t) &=& A(t)x(t)+B(t)u(t),\; x(t_0)=0, \\
y(t) &=& C(t) x(t),
 \end{array} 
 \right .
\end{equation}
where $x(t) \in \mathbb{R}^n$ is the state vector, $u(t) \in \mathbb{R}^p$  is the control and $y(t) \in \mathbb{R}^p$ is the output. The matrices $A(t) \in  \mathbb{R}^{n \times n}$, $B(t) \in  \mathbb{R}^{n \times p}$ and $C(t) \in  \mathbb{R}^{p \times n}$ are assumed to be  continuous and bounded for all $t \in [t_0,\, T_f]$. \\
The LTV dynamical system (\ref{ltv1}) can also be denoted as
\begin{equation}\label{ltv22}
{\Sigma}(t)  \equiv  
\left [
\begin{array}{c|c}
A(t) & B(t)\\
\hline
C(t) & 0
\end{array}
\right ].
\end{equation}
 In many applications, such as circuit simulation, or time dependent PDE control problems, the dimension $n$ of ${\Sigma} $  is quite large, while the number of inputs  and outputs  is small $p \ll n$. In these large-scale settings, the system dimension makes the computation infeasible due to memory, time limitations and ill-conditioning. To overcome these drawbacks, one approach consists in reducing the model. The goal is to produce a low order system that has similar response characteristics as the original system with  lower storage requirements and evaluation time.\\
 The  reduced order  dynamical system can be expressed as follows
 \begin{equation}\label{ltv2}
 {\ \Sigma}_m \left \{
 \begin{array}{lll}
{\dot x}_m(t)   =   A_m(t) x_m(t) +B_m(t)u(t) \\
 \\
 y_m(t)   =   C_m(t)x_m(t) \\
 \end{array}
 \right.
 \end{equation}
 where  $x_m  \in  \mathbb{R}^m$,  $ y_m \in  \mathbb{R}^p$, $A_m \in  \mathbb{R}^{m \times m}$, $B \in  \mathbb{R}^{m \times p}$ and $C_m \in  \mathbb{R}^{p \times m}$ with $m \ll n$. The reduced dynamical system (\ref{ltv2}) is also represented as
 \begin{equation}\label{ltv3}
 \Sigma_m(t)\equiv 
 \left [
 \begin{array}{c|c}
 A_m(t) & B_m(t)\\
 \hline
 C_m(t) & 0
 \end{array}
 \right ].
 \end{equation}
 The reduced order dynamical system should be constructed in order that
 \begin{itemize}
 \item  The  output $y_m(t)$ of the reduced system approaches  the output $y(t)$ of the original system.
 \item Some  properties of the original system such as passivity and stability   are  preserved.
 \item The computation  methods are  steady and efficient. 
 \end{itemize}

\noindent One of the well known methods for constructing such reduced-order dynamical systems is the balanced truncation method for LTV systems \cite{sandberg,shokoohi,verries}; see also \cite{glover1,moore,mullis} for the linear time-independent case.  This method requires the LTV controllability  and observability  Gramians $P(t)$ and $Q(t)$ defined as the solutions of the differential Lyapunov matrix equations
\begin{equation}
\label{control}
{\dot P}(t) = A(t)P(t)+P(t)A(t)^T +B(t)B(t)^T,\; P(t_0)=0,
\end{equation}
and
\begin{equation}
\label{obs}
{\dot Q}(t) = A^T(t)P(t)+P(t)A(t) +C(t)^TC(t),\; Q(T_f)=0.
\end{equation}
Using the formulae \eqref{solexacte1}, the differential Lyapunov equation \eqref{control} has the unique symmetric and positive solution $P(t)$ given by
$$P(t)=\int_{t_0}^t \Phi_A(t,\tau)B(\tau)B^T(\tau)\Phi^T_A(t,\tau)d\tau,  $$
where the transition matrix $\Phi_A(t,\tau)$ is the unique solution of the problem
$$\displaystyle {\dot \Phi}_A(t,\tau)=A(t) \Phi_A(t,\tau),\;\; \Phi_A(t,t)=I.$$
The observability Gramian is given by
$$Q(t)=\int_{t}^{T{_f}} \Phi^T_A(\tau,t)C^T(\tau)C(\tau)\Phi_A(\tau,t) d\tau.$$
The two  LTV controllability  and observability Gramians $P(t)$ and $Q(t)$ are then used to construct a new balanced system  such that $\tilde P(t)=\tilde Q(t)=diag(\sigma_1(t),\ldots,\sigma_n(t))$ where the Hankel singular values are given by  $\sigma_i(t)= \sqrt{\lambda_i(P(t)Q(t)}$, $i=1,\ldots,n$  and order in decreasing order.\\
The concept of balancing  is to a transform the original LTV system to an equivalent one  in which the states that are difficult to reach are also difficult to observe, which is finding an equivalent new LTV system such that the new 
 Gramians ${\widetilde P}$ and ${\widetilde Q}$ are such that 
$${\widetilde P}(t) = {\widetilde Q}(t)=diag(\sigma_1,\ldots,\sigma_n)$$
where $\sigma_i$ is the $i$-th Hankel singular value of the LTV system; i.e. $$\sigma_i = \sqrt{\lambda_i(P(t)Q(t))}.$$
Consider  the Cholesky decompositions of the Gramians $P$ and $Q$:
\begin{equation}
\label{l1}
P(t)=L_c(t)L_c(t)^T,\;\; Q(t)=L_o(t)L_o(t)^T,
\end{equation}
and consider also the singular value decomposition of $L_c(t)^T L_o(t)$ as
\begin{equation}
\label{l2}
L_c(t)^T L_o(t) = Z(t) \Sigma(t) Y(t)^T,
\end{equation}
where $Z(t)$ and $Y(t)$ are unitary $n \times n$ matrices and $\Sigma$ is a diagonal matrix containing the 
singular values.
The balanced truncation consists in determining a reduced order model  by truncating the states   corresponding to the small Hankel singular values. Under certain conditions stated in \cite{shokoohi}, one can construct the low order model ${ \Sigma}_m(t)$ as follows:  We set 
\begin{equation}
\label{red11}
V_m(t) =L_o(t) Y_m(t) \Sigma_m(t)^{-1/2}\;{\rm and }\; W_m(t) = L_c(t) Z_m(t) \Sigma_m(t)^{-1/2},
\end{equation}
where $\Sigma_m(t)=diag(\sigma_1(t),\ldots,\sigma_m(t))$; $Z_m(t)$ and $Y_m(t)$ correspond to the leading $m$ columns of the matrices $Z(t)$ and $Y(t)$ given by the singular value decomposition (\ref{l2}). The  matrices of the reduced LTV system

\begin{equation}
\label{red12}
W_m(t)^TV_m(t)A_m(t) = V_m(t)^T A(t) W_m(t)-V_m(t)^T {\dot W}_m(t),\; B_m(t)=V_m(t)^T B(t),\; C_m(t) = C(t)W_m(t).
\end{equation}
The use of Cholesky factors in the Gramians $P(t)$ and $Q(t)$ is not applicable for large-scale problems. Instead,  one can compute low rank approximations of $P(t)$ an $Q(t)$ as given by \eqref{approx} and use them to construct an approximate balanced truncation model.\\

\noindent As $A$, $B$ and $C$ are time-dependent, the direct application of the two approaches we developed is too expensive. Instead, we can apply directly   an integration method such as BDF to the differential Lyapunov matrix equations \eqref{control} and \eqref{obs}. Then, at each iteration of the BDF method, we obtain a large Lyapunov matrix equation that can be numerically solved  by using the extended block Arnoldi algorithm. \\
Consider the differential matrix equation \eqref{control}, then,  at each iteration of the BDF method, the approximation $P_{k+1}$ of  $P(t_{k+1})$ where $P$ is the exact solution of \eqref{control},  is given by the implicit relation 
\begin{equation}
\label{bdf1}
P_{k+1} = \displaystyle \sum_{i=0}^{p-1} \alpha_i P_{k-i} +h_k \beta {\mathcal G}(G_{k+1}),
\end{equation} 
where $h_k=t_{k+1}-t_k$ is the step size, $\alpha_i$ and $\beta_i$ are the coefficients of the BDF method as listed  in Table \ref{tab1} and ${\mathcal G}(X)$ is  given by 
$${\mathcal G}(X)= A^T\,X+X\,A+\,B\,B^T.$$

\noindent The approximate solution $P_{k+1}$ solves the following matrix equation
\begin{equation*}
-P_{k+1} +h_k\beta (A^T P_{k+1} + P_{k+1} A+ B B^T) + \displaystyle \sum_{i=0}^{p-1} \alpha_i P_{k-i} = 0,
\end{equation*}
which can be written as the following  continuous-time algebraic Riccati equation
\begin{equation}
\label{lyapbdf1}
\mathcal{A}_k^T\, P_{k+1}  +\,P_{k+1}\, \mathcal{A}_k  + \mathcal{B}_k\, \mathcal{B}_k^T  =0.
\end{equation}
Assuming that at each timestep, $P_k$ can be approximated as a product of  low rank factors  $P_{k}\approx \tilde Z_{k} \tilde  Z_k^T$, $\tilde Z_k \in \mathbb{R}^{n \times m_k}$, with $m_k \ll n$, the coefficient matrices are given by
$$\mathcal{A}_k= h_k\beta A -\displaystyle \frac{1}{2}I,~~\mbox{and }\mathcal{B}_{k+1}=[\sqrt{h_k\beta} B, \sqrt{\alpha_0} \tilde Z_k^T,\ldots,\sqrt{\alpha_{p-1}} \tilde Z_{k+1-p}^T]^T.$$

\noindent A good way for solving the Lyapunov matrix equation \eqref{lyapbdf1} is by  using the block or extended block Arnoldi algorithm applied to the pair $(\mathcal{A}_k,\mathcal{B}_k)$. This allows us  to obtain low rank  approximate solutions in factored forms. The procedure is as follows: applying for example the  block Arnoldi to the pair $(\mathcal{A}_k,\mathcal{B}_k)$ we get, at step $m$ of the Arnoldi process,  an orthonormal basis of the extended block Krylov subspace formed by the columns of the matrices: $\{V_{1,k},\ldots,V_{m,k}\}$ and also a block upper Hessenberg matrix $\mathbb{H}_{m,k}$. Let  $\mathbb{V}_{m,k}=[V_{1,k},\ldots,V_{m,k}]$ and $\mathbb{H}_{m,k} = \mathbb{V}_{m,k}^T \mathcal{A}_k \mathbb{V}_{m,k}$. Then the obtained low rank approximate solution to the solution $P_{k+1} $ of  \eqref{lyapbdf1} is given as $P_{m,k}=\mathbb{V}_{m,k} \mathbb{Y}_{m,k} \mathbb{V}_{m,k}^T$ where $\mathbb{Y}_{m,k}$ is solution of the following  low order Lyapunov equation 
\begin{equation}
\label{lyapred1}
\mathbb{H}_{m,k} \mathbb{Y}_{m,k}+\mathbb{Y}_{m,k} {\mathbb{H}_{m,k}}^T+ {\mathcal{\tilde B}}_k\, {{\mathcal{\tilde B}}_k}^T=0,
\end{equation}
where $\mathcal{\tilde B}_k=  \mathbb{V}_{m,k}^T \, \mathcal{B}_k$. As stated in Remark 1, the approximate solution can be given in a factored form.

\section{Numerical examples}

In this section, we compare the two approaches presented in this paper. The exponential approach (EBA-exp) summarized in Algorithm \ref{algo_EBA_exp}, which is based on the approximation of the solution to \eqref{lyap1} applying a quadrature method to compute the projected exponential form solution \eqref{gm}.  We used a scaling and squaring strategy, implemented in the MATLAB  \textbf{expm} function; see  \cite{higham05,moler03} for more details.  The second method (Algorithm \ref{algo_EBA_BDF}) is based on the BDF  integration method applied to the projected Lyapunov equation \eqref{lowlyap}. The basis of the projection subspaces were generated by the extended block Arnoldi algorithm for both methods.
All the experiments were performed on a laptop with an  Intel Core i7 processor and 8GB of RAM. The algorithms were coded in Matlab R2014b. \\

\noindent {\bf Example 1}. 
The matrix $A$  was   obtained from the 5-point discretization of the operators 
\begin{equation*}
L_A=\Delta u-f_1(x,y)\frac{\partial u}{\partial x}+ f_2(x,y)\frac{\partial u}{\partial y}+g_1(x,y),
\end{equation*}
on the unit square $[0,1]\times [0,1]$ with homogeneous Dirichlet boundary conditions.  The number of inner grid points in each direction is  $n_0=$ and the dimension of the matrix $A$ was $n = n_0^2=$. Here we set $f_1(x,y) = 10xy$, $f_2(x,y)= e^{x^2y}$, $f_3(x,y) = 100y$, $f_4(x,y)= {x^2y}$ ,  $g_1(x,y) = 20y$ and $g_2(x,y)=x\,y$.  The time interval considered was $[0,\,2]$ and the initial condition $X_0=X(0)$ was choosen as the low rank product $X_0=Z_0Z_0^T$, where $Z_0=0_{n \times 2}$.
For both methods, we used  projections onto the Extended Block Krylov subspaces $${\mathcal K}_k(A,B) = {\rm Range}(B,A\,B,\ldots,A^{m-1}\,B,A^{-1}\,B,\ldots,(A^{-1})^m\,B)$$ and the  tolerance was set to $10^{-10}$ for the stop test on the residual.  For the EBA-BDF method,  we used a 2-step BDF scheme with a constant timestep $h$. The entries of the matrix $B$ were random values uniformly distributed on the interval $[0, \, 1]$ and the number of the columns in $B$  was $s=2$. \\
 literature. In order to check if our approaches produce reliable results, we began comparing our results to the one given by Matlab's ode23s solver which is designed for stiff differential equations. This was done by vectorizing our DLE, stacking the columns of $X$ one on top of each other. This method, based  on Rosenbrock integration scheme, is not suited to large-scale problems. Due to the memory limitation of our computer when running the ode23s routine,  we chose a  size of $100\times 100$ for the matrix $A$.\\

\begin{figure}[H]
	\begin{center}
		\includegraphics[width=15cm,height=7cm]{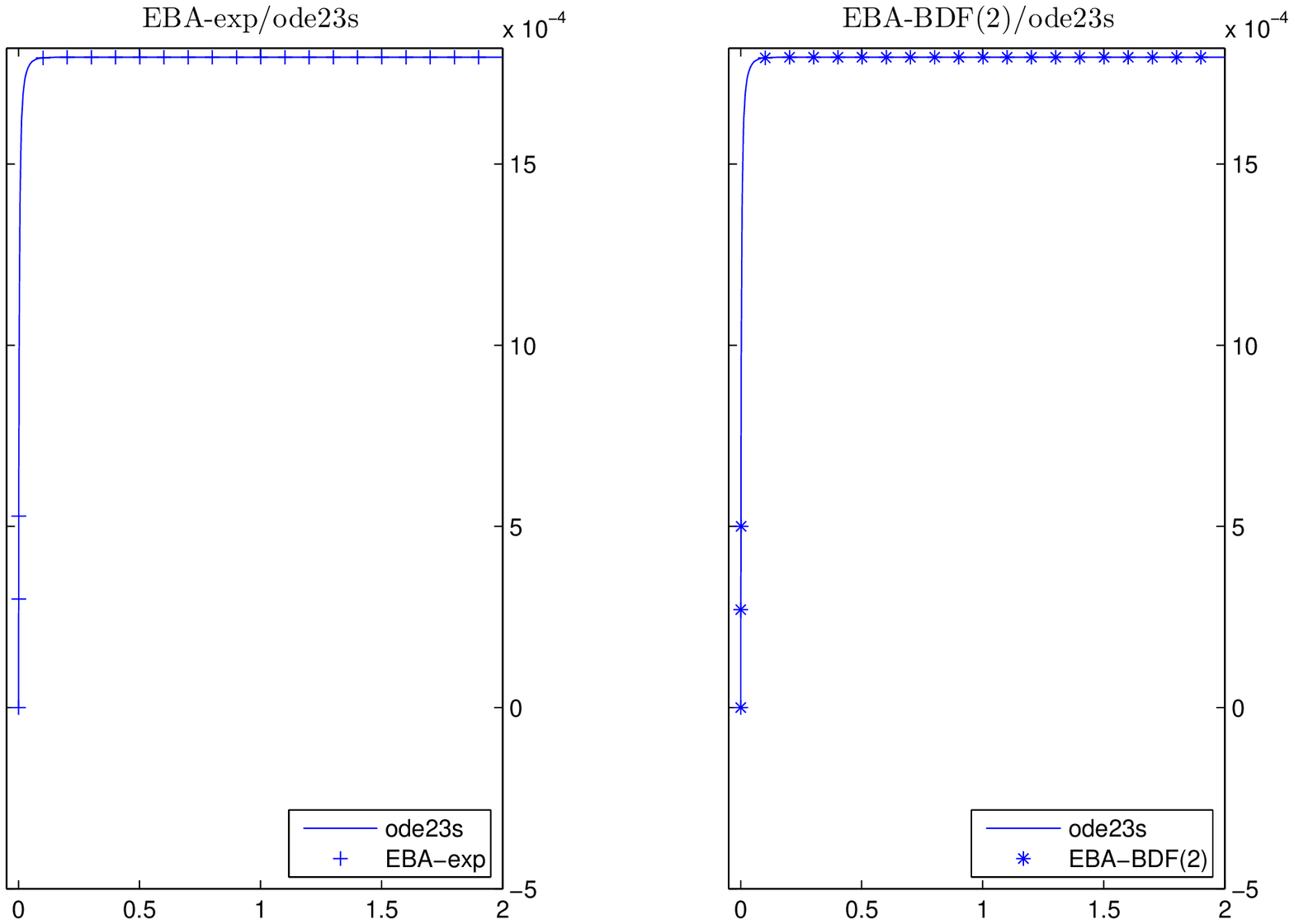}
		\caption{ Values of $X_{11}(t)$ for  $t \in [0,\, 2]$}\label{Figure1}
	\end{center}
\end{figure}

\noindent In Figure \ref{Figure1}, we compared the component $X_{11}$ of the solution obtained by the methods tested in this section, to the solution provided by the ode23s method from Matlab,  on the time interval $[0,\,2]$, for $size(A)=100\times 100$ and a constant timestep $h=10^{-3}$.
\noindent We observe that all the considered methods give similar results in terms of accuracy. The relative error norms   $\displaystyle{\frac{\|X_{EBA-exp}(t_f)-X_{ode23s}(t_f)\|}{\|X_{ode23s(t_f)}\|}}$ and $\displaystyle{\frac{\|X_{EBA-BDF(2)}(t_f)-X_{ode23s}(t_f)\|}{\|X_{ode23s(t_f)}\|}}$ at final time $t_f=2$ were equal to $1.8\times 10^{-10}$ and  $9.1\times 10^{-11}$ respectively.  The runtimes were respectively 0.59s, 5.1s for the EBA-exp and EBA-BDF(2) methods and 1001s for the ode23s routine.

\noindent In Table \ref{tab2}, we give  the obtained runtimes in seconds, for the resolution of Equation \eqref{lyap1}  for $t \in [0,\, 2]$, with a timestep $h=0.001$ and the Frobenius norm of the residual at the final time.

\begin{table}[h!]
	\begin{center}
		\begin{tabular}{c | c c c }
			size($A$)&EBA-exp&EBA-BDF(2)&Residual norm\\
			\hline
			$2500 \times 2500$&$3.23$ s&$31.9$ s&$\mathcal{O}(10^{-9})$~$(m=16)$\\
			$6400 \times 6400$&$5.2$ s&$81.6$ s&$\mathcal{O}(10^{-9})$~$(m=19)$\\
			$10000 \times 10000$&$5.6$ s&$168$ s&$\mathcal{O}(10^{-8})$~$(m=19)$\\
			$22500 \times 22500$&$11.8$ s&$1546$ s&$\mathcal{O}(10^{-8})$~$(m=23)$\\
			\hline
		\end{tabular}
		\caption{runtimes and residual norms for EBA-exp and EBA+BDF(2)}\label{tab2}
	\end{center}
\end{table}

\begin{figure}[H]
	\begin{center}
		\includegraphics[width=7cm,height=5cm]{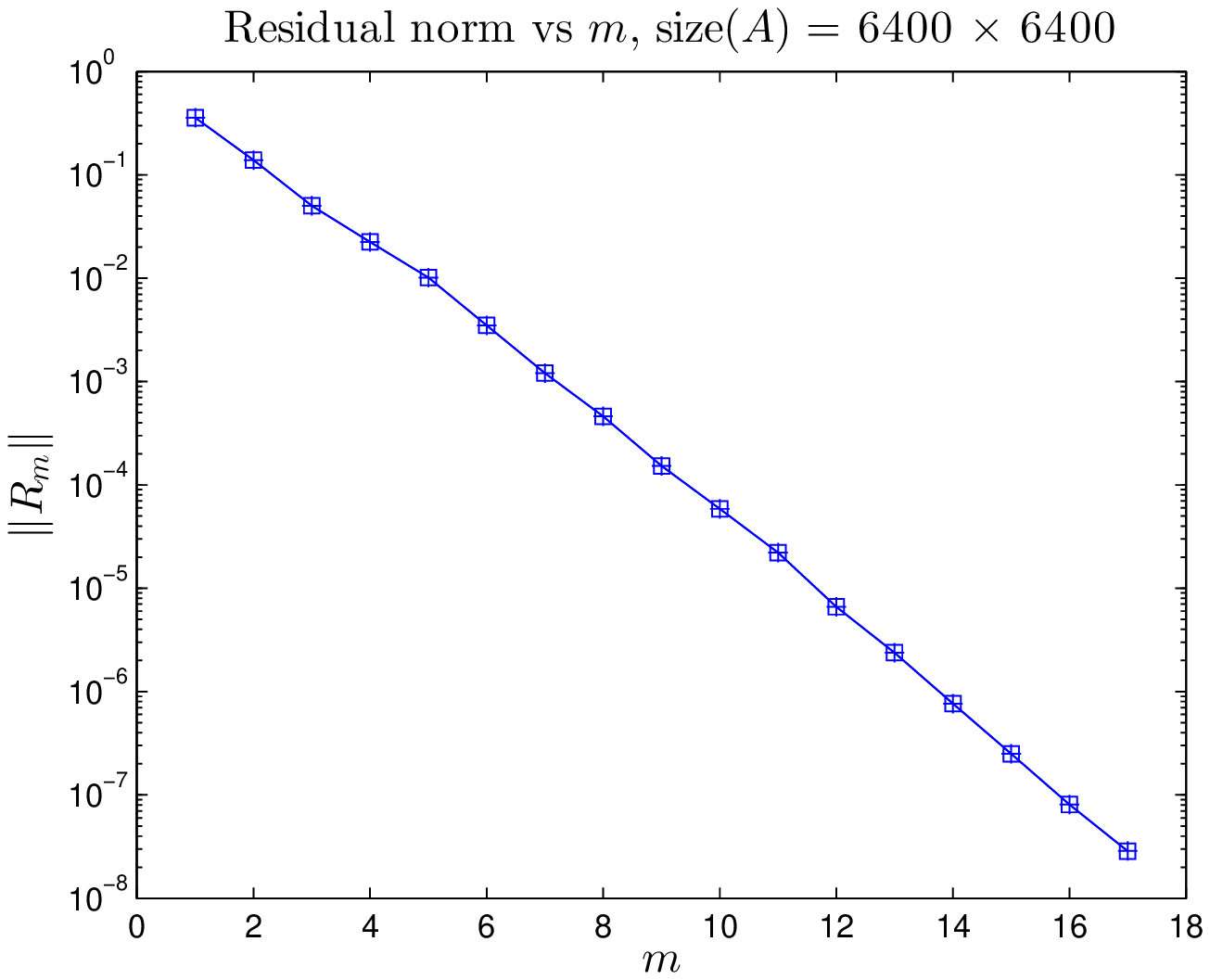}
		\caption{ Residual norms \textit{vs} the number of extended Arnoldi iterations $m$}\label{Figure2}
	\end{center}
\end{figure}

\noindent The results  in Table \ref{tab2} illustrate that the EBA-exp method clearly outperforms the EBA-BDF(2) method in terms of computation time even though both methods are equally accurate. In Figure \ref{Figure2}, we featured the norm of the residual at final time $t=2$ for both EBA-exp and EBA-BDF(2) methods  for size($A$)$=6400 \times 6400$ in function of the number $m$ of extended Arnoldi iterations. We observe that the plots coincide for both methods. 

\noindent {\bf Example 2}. 
This example  comes from the autonomous linear-quadratic
optimal control problem of one dimensional heat flow
\begin{eqnarray*}
	\frac{\partial}{\partial t} x(t,\eta) & = & \frac{\partial^2}{\partial \eta^2} x(t,\eta)+b(\eta) u(t)   \\
	x(t,0) & = & x(t,1)=0, t>0\\
	x(0,\eta) & = & x_0(\eta), \eta \in [0,1]\\
	y(x) & = & \int_0^1 c(\eta) x(t,\eta) d \eta, x>0.
\end{eqnarray*}
Using a standard finite element approach based on the first order B-splines, we obtain the following ordinary differential equation
\begin{eqnarray}\label{ode1}
M \dot {\tt x}(t)& = & K {\tt x}(t) + F u(t)\\
y(t) & =& C {\tt x}(t),
\end{eqnarray}
where the matrices $M$ and $K$ are given by:
$$M=\frac{1}{6n}\left(
\begin{array}{ccccc}
4&1&&&\\
1&4&1&&\\
&\ddots&\ddots&\ddots\\
&&1&4&1\\
&&&1&4
\end{array}
\right), \;\; K=-\alpha\,n\,\left(
\begin{array}{ccccc}
2&-1&&&\\
-1&2&-1&&\\
&\ddots&\ddots&\ddots\\
&&-1&2&-1\\
&&&-1&2
\end{array}
\right).$$
Using the semi-implicit Euler method, we get  the following discrete dynamical system
$$(M-\Delta t K)\, \dot{x}(t)= M\, x(t)+ \Delta t\, F u_k.$$
We set $A=(M-\Delta t K)^{-1}\, M$ and $B=\Delta t\, (M-\Delta t K)^{-1}\, F$. The entries of the   $ n \times s$  matrix $ F $ and the $s \times n$ matrix  $ C $ were random values uniformly distributed on $ [0,\,1] $. In our experiments we used  $n=$, $s=2$,   $\Delta t = 0.01$   and  $\alpha=0.05$.

\noindent In Table \ref{tab3}, we give  the obtained runtimes in seconds, for the resolution of Equation \eqref{lyap1}  for $t \in [0,\, 2]$, with a timestep $h=0.001$ and the Frobenius norm of the residual at the final time.

\begin{table}[h!]
	\begin{center}
		\begin{tabular}{c | c c c }
			size($A$)&EBA-exp&EBA-BDF(2)&Residual norms\\
			\hline
			$2500 \times 2500$&$1.0$ s&$8.0$ s&$\mathcal{O}(10^{-11})$~$(m=11)$\\
			$6400 \times 6400$&$4.9$ s&$14.4$ s&$\mathcal{O}(10^{-14})$~$(m=11)$\\
			$10000 \times 10000$&$11.5$ s&$29.7$ s&$\mathcal{O}(10^{-13})$~$(m=11)$\\
			$20000 \times 20000$&$11.8$ s&$173.4$ s&$\mathcal{O}(10^{-13})$~$(m=11)$\\
			\hline
		\end{tabular}
		\caption{runtimes and residual norms for EBA-exp and EBA-BDF(2)}\label{tab3}
	\end{center}
\end{table}

\begin{figure}[H]
	\begin{center}
		\includegraphics[width=7cm,height=5cm]{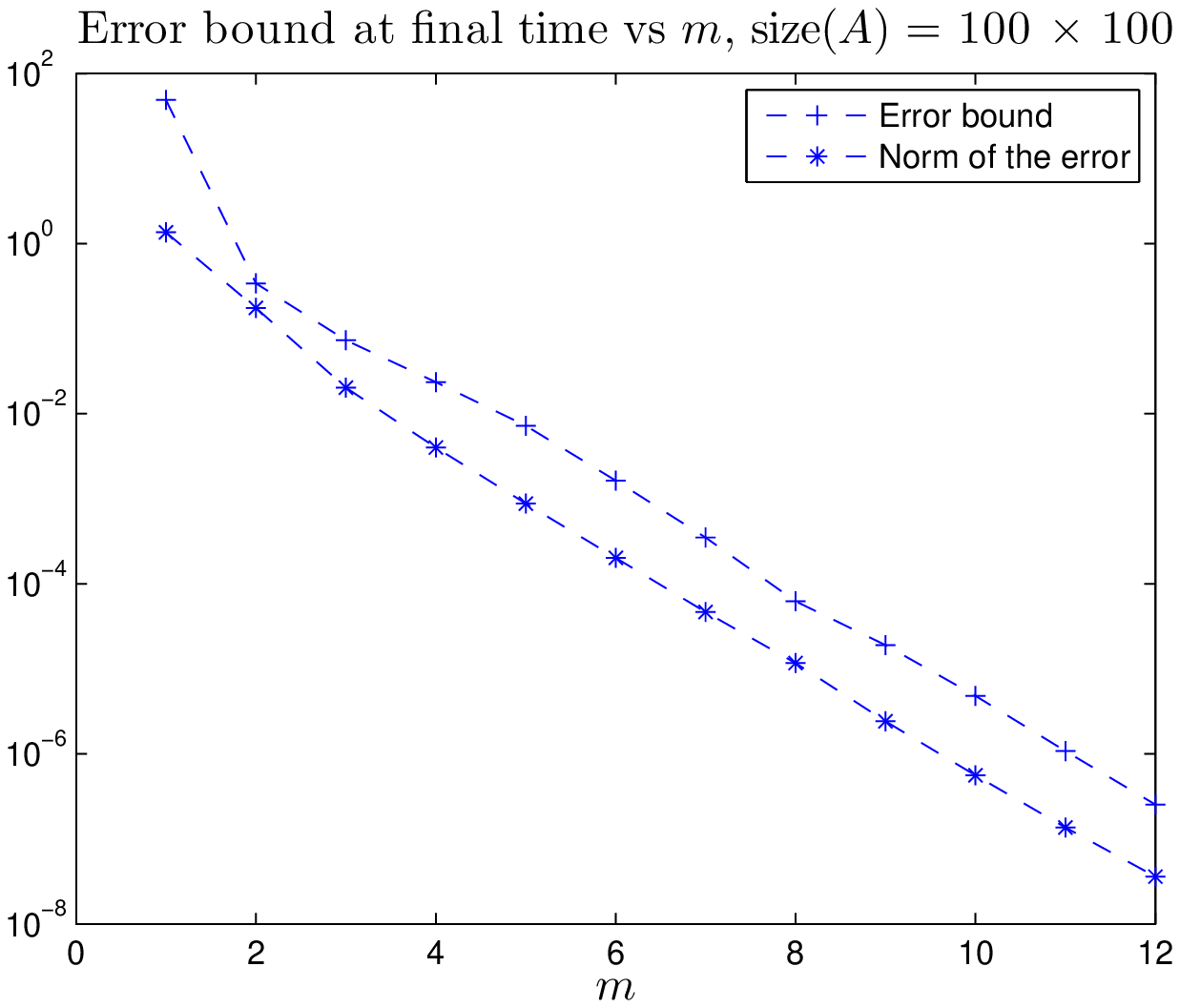}
		\caption{ Upper bounds of the error norms and   computed error norms \textit{vs}  the number of iterations}\label{Figure3}
	\end{center}
\end{figure}

\noindent The figures in Table \ref{tab3} illustrate the gain of speed provided by the EBA-exp method. Again, both methods performed similarly in terms of accuracy. In figure \ref{Figure3}, we considered the case size($A$)$=100\times 100$ and  plotted the upper bound of the error norms as stated in Formula \eqref{upperbound} at the final time $T_f$ against the computed norm of the errors, taking the solution given by the integral formula \eqref{solexacte3} as a reference, in function of the number $m$ of Arnoldi iterations for the EBA-exp method. 

\noindent {\bf Example 3} 
In this last example, we applied the EBA-BDF(1) method to the well-known problem Optimal Cooling of Steel Profiles. The matrices were extracted from the IMTEK collection \footnote{https://portal.uni-freiburg.de/imteksimulation/downloads/benchmark}. We compared the EBA-BDF(2) method to the EBA-exp method for problem sizes $n=1357$ and $n=5177$, on the time interval $[0\,,1000]$. The initial value $X_0$ was chosen as $X_0=0$ and  the timestep was set to $h=0.01$. The tolerance for the Arnoldi stop test was set to $10^{-7}$ for both methods and the projected low dimensional Lyapunov equations were numerically solved by the  solver ({\tt lyap} from Matlab) at each iteration of the extended block Arnoldi  algorithm for the EBA-BDF(2) method.
\vspace{0.2cm}

\begin{table}[h!]
	\begin{center}
		\begin{tabular}{c | c c c }
			size($A$)&EBA-exp&EBA-BDF(2)&Residual norms\\
			\hline
			$1357\times1357$&$11.1$ s&$515.5$ s&$\mathcal{O}(10^{-8})$~$(m=6)$\\
			$5177 \times 5177$&$148.7$ s&$1721$ s&$\mathcal{O}(10^{-7})$~$(m=39)$\\
			\hline
		\end{tabular}
		\caption{Optimal Cooling of Steel Profiles: runtimes and residual norms for EBA-exp and EBA-BDF(2)}\label{tab4}
	\end{center}
\end{table}
\noindent In Table \ref{tab4},  we listed the obtained runtimes which again showed the advantage of the EBA-exp method in terms of execution time and similar accuracy for both methods.

\section{Conclusion}
We presented in the present paper two new approaches for computing approximate solutions to large scale differential  Lyapunov matrix equations. The first one comes naturally from the exponential expression of the exact solution and the use of approximation techniques of the exponential of a matrix times a block of vectors. The second approach is obtained by first projecting the initial problem onto a block Krylov (or extended Krylov) subspace, obtain a low dimensional differential Lyapunov equation which is solved by using the well known BDF integration method. We gave some theoretical results such as the exact expression of the residual norm and also upper bounds for the norm of the errors. An application in model reduction for  linear time-varying dynamical  systems is also given. Numerical experiments show that both methods are promising for large-scale problems, with a clear advantage for the EBA-exp method in terms of computation time.

\bibliographystyle{plain}
\end{document}